\documentclass[12pt]{amsart}
\usepackage{amssymb}

\newtheorem{thm}{Theorem}[section]
\newtheorem{cor}[thm]{Corollary}
\newtheorem{prop}[thm]{Proposition}
\newtheorem{lemma}[thm]{Lemma}

\newtheorem{conj}[thm]{Conjecture}

\theoremstyle{remark}
\newtheorem{remark}[thm]{Remark}
\newtheorem{example}[thm]{Example}

\theoremstyle{definition}
\newtheorem{defn}[thm]{Definition}

\numberwithin{equation}{section}

\newcommand{\bbA}{{\mathbb A}}
\newcommand{\bbC}{{\mathbb C}}

\newcommand{\bbL}{{\mathbb L}}
\newcommand{\bbR}{{\mathbb R}}

\newcommand{\bbZ}{{\mathbb Z}}

\newcommand{\bbP}{{\mathbb P}}

\newcommand{\cF}{{\mathcal F}}

\newcommand{\cO}{{\mathcal O}}
\newcommand{\cP}{{\mathcal P}}

\newcommand{\cM}{{\mathcal M}}

\newcommand{\cD}{{\mathcal D}}
\newcommand{\cV}{{\mathcal V}}
\newcommand{\cA}{{\mathcal A}}
\newcommand{\cB}{{\mathcal B}}

\newcommand{\cC}{{\mathcal C}}
\newcommand{\cE}{{\mathcal E}}

\newcommand{\cU}{{\mathcal U}}

\newcommand{\cT}{{\mathcal T}}

\newcommand{\D}{\operatorname{D}}

\newcommand{\Hom}{\operatorname{Hom}}
\newcommand{\Ext}{\operatorname{Ext}}

\newcommand{\shHom}{\underline{\operatorname{Hom}}}

\newcommand{\add}{\operatorname{add}}
\newcommand{\smd}{\operatorname{smd}}

\newcommand{\Ind}{\operatorname{Ind}}
\newcommand{\Res}{\operatorname{Res}}

\newenvironment{pf}{\smallskip\noindent{\bf Proof.}\ }{\qed\smallskip}

\title{Grothendieck ring of pretriangulated categories}
\author{Alexey I.~Bondal}
\address{Steklov Mathematical Institute, 8 Gubkina St. Moscow, Russia}
\email{bondal@mi.ras.ru}
\author{Michael Larsen}
\address{Department of Mathematics, Indiana University,
Bloomington, IN 47405, USA}
\email{larsen@math.indiana.edu}
\author{Valery A.~Lunts}
\address{Department of Mathematics, Indiana University,
Bloomington, IN 47405, USA}
\email{vlunts@indiana.edu}

\thanks{The first named author was partially supported by the
CRDF grant RM1-2405-MO-02.
The second named author was partially supported by NSF grant DMS-0100537.
The third named author was partially supported by NSA grant MDA904-01-1-0020
and CRDF grant RM1-2405-MO-02}

\dedicatory{To the blessed memory of Andrej Nikolaevich Tyurin}

\begin{document}

\begin{abstract} We consider the abelian group $\cP\cT$ generated by quasi-equivalence classes
of pretriangulated DG categories with relations coming from semi-orthogonal decompositions
of corresponding triangulated categories.
We introduce an operation of "multiplication" $\bullet$ on the collection of DG categories
which makes this abelian group into a commutative ring. A few applications are considered:
representability of "standard" functors between derived categories of coherent sheaves
on smooth projective varieties and a construction of an interesting motivic measure.
\end{abstract}

\maketitle

\section*{Introduction}

This work grew out of an attempt to construct the Grothendieck ring of (equivalence classes
of) triangulated categories. Namely, one considers an abelian group $\cT$
generated by equivalence classes of triangulated categories with relations coming from
semiorthogonal decompositions: $[\cA]=[\cB]+[\cC]$ if $\cA$ admits a semiorthogonal
decomposition with the two summands being equivalent to $\cB$ and $\cC$ respectively.
We wanted to define a product $\bullet$ of triangulated categories which would make $\cT$ a
commutative associative ring. This answer was supposed to be known in the
following situations:

1) If $A$ and $B$ are nice algebras (say, finite dimensional and of finite cohomological
dimension), and $D(A)$ and $D(B)$ are the bounded derived categories of (finite) $A$- and
$B$-modules respectively, then $D(A)\bullet D(B)=D(A\otimes B)$.

2) If $X$ and $Y$ are smooth projective varieties and $D(X)$ and $D(Y)$ are the corresponding
bounded derived categories of coherent sheaves, then $D(X)\bullet D(Y)=D(X\times Y)$.

After a while we came to conclusion that triangulated categories are not rigid enough
(a well known problem being that taking the cone of a morphism is not a functorial operation),
and convinced ourselves to work with {\it pretriangulated} categories. These
are DG categories where the cone of a morphism {\it is} a functor, so they are more rigid and
easier to work with than the triangulated categories. In a sense this is going back from
the homotopy category to the abelian category of complexes. The {\it homotopy} category
of a pretriangulated category is a triangulated category. So this approach gives a
satisfactory solution to our original problem.

For example, if $X$ and $Y$ are smooth projective varieties and $I(X)$ and $I(Y)$ are
the DG categories of bounded below injective complexes of $\cO _X$- and $\cO _Y$-modules
respectively with
bounded coherent cohomology, then our product gives
$$I(X)\bullet I(Y)=I(X\times Y).$$

We present two applications of this last formula. The first application is the
representability of {\it standard} functors $F:D(X)\to D(Y)$. Namely, we prove that if
an exact functor $F$ comes from a DG functor between standard enhancements of $D(X)$
and $D(Y)$ respectively, then there is an object $P\in D(X\times Y)$ which represents
$F$: there exists an isomorphism of functors
$$F(\cdot )=\bbR q_*(p^*(\cdot)\stackrel{\bbL}{\otimes}P),$$
where $X\stackrel{p}{\leftarrow}X\times Y\stackrel{q}{\rightarrow}Y$ are the projections.
We also conjecture that every exact functor between $D(X)$ and $D(Y)$ is standard.

Our second application is the construction of an interesting {\it motivic measure}, i.e.
a homomorphism from
the Grothendieck ring of varieties to the Grothendieck ring of (quasi-equivalence classes)
of pretriangulated categories, by sending a smooth projective variety $X$ to $I(X)$.
We show that the kernel of this homomorphism contains the element $\bbL-1$, where $\bbL$
is the class of the affine line.

The authors would like to thank Dima Orlov, Bernhard Keller and Vladimir Drinfeld for
useful discussions and suggestions.

\section{Generators and semiorthogonal decompositions of triangulated  categories}

Fix a field $k$. All categories and all functors that we will consider are
assumed to be $k$-linear. Denote by $Vect$ the category of $k$-vector spaces.

For a smooth projective variety $X$ over $k$ we denote $D(X)=D^b(coh_X)$ the bounded
derived category of coherent sheaves on $X$.

In this section we recall some results about triangulated categories following
[Bo],[BoKa1],[BoVdB],[BN],[KeVo]. In the end we define the Grothendieck group $\cT$ of
triangulated categories.

\subsection{Generators and representability of cohomological functors}

Let $\cE =(E_i)_{i\in I}$ be a class of objects in a triangulated category $\cA$.
A {\it triangulated envelope} of $\cE$ is the smallest strictly full triangulated
subcategory of $\cA$ which contains $\cE$.

If $\cA$ is a triangulated category then a triangulated subcategory $\cB \subset \cA$ is
called epaisse (thick) if it is closed under isomorphisms and direct summands.

We say that $\cE$
{\it classically generates} $\cA$ if the smallest {\it epaisse triangulated} subcategory of
$\cA$ containing $\cE$ (called the {\it epaisse envelope} of $\cE$ in $\cA$) is equal to
$\cA$ itself. We say that $\cA$ is {\it finitely generated} if it is classically generated
by one object.

By the {\it right orthogonal} $\cE ^\perp$ in $\cA$ we denote the full subcategory of $\cA$
whose objects $A$ have the property $\Hom (E_i[n],A)=0$ for all $i$ and $n$. Similarly, we
define the left orthogonal ${}^\perp\cE$. Clearly $\cE^\perp$ (and ${}^\perp\cE$) is an
epaisse subcategory of $\cA$. We say that $\cE$ {\it generates} $\cA$ if $\cE ^\perp=0$.
Clearly if $\cE$ classically generates $\cA$ then it generates $\cA$, but the
converse is false.

Denote by $\add (\cE)$ the minimal strictly full subcategory of $\cA$ which contains $\cE$
and is closed under taking finite direct sums and shifts. Also denote by $\smd(\cE)$
the minimal strictly full subcategory which contains $\cE$ and is closed under taking
(possible) direct summands.

There exists a natural multiplication on the set of strictly full
subcategories of $\cA$. If $\cC$ and $\cD$ are two such
subcategories, let $\cC \star \cD$ be the strictly full
subcategory whose objects $S$ occur in exact triangles $C\to S\to D$ with
$C\in \cC$, $\D\in \cD$. This multiplication is associative in view of the octahedral axiom.
If $\cC$ and $\cD$  are closed under direct sums and/or shifts, then so is $\cC \star \cD$.

Now we define a new multiplication operation on the set of strictly full subcategories
which are closed under finite direct sums by the formula:
$$\cC \diamond \cD=\smd (\cC \star \cD).$$
One can check that this operation is associative. Put
$$\langle \cE \rangle _s=\add (\cE)\diamond ...\diamond \add (\cE) \ \ \ (\text{$s$ factors}).$$
$$\langle \cE \rangle =\bigcup_s\langle \cE \rangle _s$$
Thus $\langle \cE \rangle $ is the epaisse envelope of $\cE $ in $\cA$. So $\cE$ classically
generates $\cA$ if and only if $\langle \cE \rangle =\cA$.

\begin{defn}  $\cE$ {\it strongly generates} $\cA$ if $\cA =\langle \cE \rangle _s$
for some $s$. We say that $\cA$ is {\it strongly finitely generated} if it is strongly
generated by one object.
\end{defn}

\begin{remark} Let $F:\cA \to \cB$ be an exact functor between triangulated categories,
which is surjective on isomorphism classes of objects. Assume that $\cA$ is classically
(resp. strongly) generated by a collection $\cE \subset \cA$. Then $\cB$ is classically
(resp. strongly) generated by the collection $F(\cE)$.
\end{remark}

\begin{remark} If a triangulated category $\cA$ has a strong generator, then any classical
generator of $\cA$ is a strong one.
\end{remark}

 A triangulated category $\cA$ is called
$\Ext$-finite if $\dim \oplus _n\Hom _{\cA}(A,B[n])<\infty$ for any
two objects $A$ and $B$.

\begin{defn} Let $\cA$ be an $\Ext$-finite triangulated category. A  functor $h:\cA
\to Vect$ is called cohomological if it takes exact triangles to
long exact sequences. We say that $h$ is of finite type if for
every object $B\in \cA$ the vector space $\oplus_nh(B[n])$ is finite dimensional.
Each object $A\in \cA$ gives rise to a
cohomological functor $h_A(\cdot):=\Hom(A,\cdot)$. Similarly,
$h^A(\cdot):=\Hom (\cdot ,A)$ defines a contravariant
cohomological functor. Cohomological functors isomorphic to $h_A$
or $h^A$ are called representable. Note that representable
functors are of finite type. The category $\cA$ is called left
(resp. right) saturated if every cohomological (resp.
contravariant cohomological) functor of finite type $h:\cA \to
Vect$ is representable. We call $\cA$ saturated if it is left and
right saturated.
\end{defn}

Recall that a category is called Karoubian, if every projector splits.

There is a simple criterion for a triangulated category to be saturated.

\begin{thm} [BoVdB] Assume that a triangulated category $\cA$ is $\Ext$-finite,
has a strong generator and is Karoubian. Then $\cA$ is right saturated.
\end{thm}

For a smooth algebraic variety $X$ over $k$ denote by $D(X)$ the bounded derived category
of coherent sheaves on $X$.

\begin{thm} [BoVdB] Let $X$ be a smooth  variety over $k$. Then $D(X)$ is
Karoubian and has a strong generator.
\end{thm}

\begin{cor} Let $X$ be a smooth projective variety over $k$. Then $D(X)$ is right saturated.
\end{cor}

\begin{remark} For a smooth projective $X$ the category $D(X)$ is equivalent to its opposite:
the functor
$${\bf D} (\cdot )=\bbR \shHom (\cdot ,\cO _X)$$
is an anti-involution of the category $D(X)$. It follows that $D(X)$ is
also left saturated.
\end{remark}

Let us also mention a few results which will be used later.

\begin{thm} [BoVdB] Let $X$ be a smooth projective variety over $k$. If an object
$E\in D(X)$ generates $D(X)$, then it strongly generates $D(X)$.
\end{thm}

\begin{thm} [BoVdB] Let $X$, $Y$ be smooth projective varieties, and let  $E^X\in D(X)$,
$E^Y\in D(Y)$ be respective generators. Then $E^X\boxtimes E^Y$ is a generator of
 $D(X\times Y)$.
 \end{thm}

\begin{thm} [BN] Assume that $\cA$ contains countable direct sums. Then the
category $\cA$ is Karoubian.
\end{thm}

\subsection{Serre functors}

Let $\cA$ be an $\Ext$-finite triangulated category. Recall [BoKa1] that a covariant
auto-equivalence  $S:\cA\to \cA$ is called a Serre functor if there exists
an isomorphism of bi-functors
$$\Hom(A,B)=\Hom(B,S(A))^*,$$
for $A,B\in \cA$.
If the Serre functor exists it is unique up to an isomorphism and is an exact functor.

\begin{example} If $X$ is a smooth projective variety of dimension $n$, the Serre functor
on the category $D(X)$ is $S(\cdot)=S_{X}(\cdot)=(\cdot ) \otimes \omega _X[n]$, where $\omega _X$ is the canonical line bundle
on $X$.
\end{example}

\begin{remark} Let $F:\cA \to \cB$, $G:\cB \to \cA$ be functors between two
$\Ext$-finite triangulated
categories with the Serre functors $S_{\cA}$ and $S_{\cB}$ respectively. Assume that the
the functor $F$ is the right adjoint to $G$, then the functor $S_{\cB}^{-1}FS_{\cA}$ is
the left adjoint to $G$.
\end{remark}

\subsection{Semiorthogonal decompositions}

Let us recall some definitions and results from [BoKa1] and [Bo].

\begin{defn} Let $\cA$ be a triangulated category, $\cB \subset \cA$ --
a strictly full triangulated subcategory. We call $\cB$ {\it right admissible}
(resp. {\it left admissible}) if for every $A\in \cA$ there exists an exact triangle
$A_{\cB}\to A \to A_{\cB ^\perp}$ (resp. $A_{{}^\perp\cB}\to A \to A_{\cB}$)
with $A_{\cB}\in \cB$ and
$A_{\cB ^\perp}\in \cB^\perp$ (resp. $A_{{}^\perp \cB }\in {}^\perp\cB$).
A subcategory is called {\it admissible} if it is both left and right
admissible.
\end{defn}

Clearly a strictly full triangulated subcategory $\cB \subset \cA$ is right (resp. left)
 admissible if and only if $\cB ^\perp$ (resp. ${}^\perp\cB$) is left (resp. right) admissible.

\begin{lemma} Let $\cB \subset \cA$ be right (resp. left) admissible. Then
${}^\perp (\cB ^\perp)=\cB$ (resp. $({}^\perp \cB) ^\perp=\cB$).
\end{lemma}

The next proposition appeared in [KeVo],[Bo],[BoKa].

\begin{prop} Let $\cA$ be a triangulated category, $\cB \subset \cA$ -- a strictly
full triangulated subcategory.  The following conditions are equivalent:

a) $\cB$ is right (resp. left ) admissible in $\cA$;

b) the embedding functor $i:\cB \to \cA$ has a right (resp. left) adjoint $i^!$ (resp. $i^*$).

It these hold, then the compositions $i^!\cdot i$ and $i^*\cdot i$ are isomorphic to
the identity functor on $\cB$.
\end{prop}

\begin{pf} (A sketch. See [KeVo] or [Bo] for details).  If $\cB$ is right admissible then
for each $A\in \cA$ the triangle
$A_{\cB}\to A \to A_{\cB ^\perp}$ is unique up to an isomorphism. Moreover, the correspondence
$A\mapsto A_{\cB}$ extends to an exact functor from $\cA$ to $\cB$. This functor is $i^!$.
Conversely, given a left adjoint $i^!$ to $i$ for any $A\in \cA$ consider the adjunction
morphism $\alpha _A:i\cdot i^!A\to A$. The required exact triangle is
$$i\cdot i^!A\to A\to C(\alpha _A).$$
Similarly for left admissible $\cB$.
\end{pf}

\begin{cor} a) Let $\cC $ be a right admissible subcategory of $\cB$, which is a right
admissible subcategory of $\cA$. Then $\cC$ is right admissible in $\cA$. Similarly for
left admissible subcategories. b) Let $\cC $ be a right admissible subcategory of $\cA$.
Assume that $\cB$ is a full subcategory of $\cA$ which contains $\cC$. Then $\cC$ is right
admissible in $\cB$.
\end{cor}

\begin{pf} a) Indeed, let $i:\cC \to \cB$, $j:\cB \to \cA$ be the embeddings and $i^!$, $j^!$ be
the corresponding right adjoints. Then $(j\cdot i)^!=i^!\cdot j^!$. b) Indeed, the
restriction of the adjoint functor of the embedding functor $\cC\to \cA$ to the subcategory
$\cB$ is the adjoint to the embedding $\cC\to \cB$.
\end{pf}

\begin{remark} Let $i:\cB\to \cA$ be an inclusion of a right (resp. left) admissible subcategory.
Assume that $\cE \subset \cA$ (strongly) generates $\cA$. Then $i^!(\cE)$ (resp. $i^*(\cE))$)
(strongly) generates $\cB$. Indeed, this follows from Remark 1.2 since the functors $i^!$ and
$i^*$ are surjective.
\end{remark}

\begin{remark} A right or left admissible subcategory $\cB \subset \cA$ is epaisse.
Hence if $\cA$ is Karoubian, then $\cB$ is also such.
\end{remark}

\begin{lemma} Let $\cA$ be a triangulated category, $\cB \subset \cA$ -- an epaisse
 triangulated
subcategory. Put $\cC =\cB ^\perp$. Let $\cA$ be classically
generated by a collection of objects $\cE \subset \cA$. Assume
that for each $E\in \cE$ there is an exact triangle $$E_{\cB}\to
E\to E_{\cC},$$
where $E_{\cB }\in \cB$ and $E_{\cC}\in \cC$.
Then the subcategory $\cB$ is right admissible
in $\cA$. Similarly for left admissible categories. Moreover, if $\cD \subset \cC$ is the
strictly full subcategory consisting of all objects $A_{\cC}$ for $A\in \cA$, then $\cD =\cC$.
\end{lemma}

\begin{pf} The proof is similar to the proof of Proposition 1.16 above.
Namely, the association $E\mapsto E_{\cB}$ can be extended to a functor from the
full subcategory, whose objects is the collection $\cE$, to $\cB$.
This functor extends to an exact functor from the triangulated envelope
of $\cE$ to $\cB$. Finally, since $\cB$ is epaisse, this functor extends to the epaisse
envelope of $\cE$, i.e. to $\cA$. This functor $\cA \to \cB$ is the right adjoint to the
inclusion functor $\cB \to \cA$. So $\cB$ is right admissible in $\cA$ by Proposition 1.16

Let us prove the last assertion of the lemma. Let $C\in \cC$ and consider the canonical
triangle
$$C_{\cB}\to C\to C_{\cC}.$$
The map $C_{\cB}\to C$ is zero, hence $C_{\cC}$ is isomorphic to $C\oplus C_{\cB}[1]$, which
implies that $C_{\cB}=0$. Hence $C\in \cD$.
\end{pf}

The following two statements establish a relation between saturated and admissible categories.

\begin{lemma}[BoKa1] Let $\cA$ be a right (resp. left) saturated triangulated category, and
$\cB \subset \cA$ be a left (resp. right) admissible subcategory. Then $\cB$ is right
(resp. left) saturated.
\end{lemma}

\begin{prop} [BoKa1] Let $\cB$ be a strictly full triangulated subcategory of a
triangulated category $\cA$ of finite type. Assume that $\cB$ is right (resp. left) saturated.
Then $\cB$ is right (resp. left) admissible in $\cA$.
\end{prop}

\begin{defn} Let $\cA$ be a triangulated category. We say that admissible subcategories
$\cB_1,...,\cB_n\subset \cA$ form a {\it semiorthogonal decomposition} of
$\cA$, denoted $\cA =(\cB _1,...,\cB _n)$ if $\cA$ is the triangulated envelope of
$\cup_iOb\cB _i$,
and $\cB _j\subset {}^\perp \cB _i$ for $j>i$.
We call $\cB _i$'s {\it semiorthogonal summands} of $\cA$.
\end{defn}

\begin{example} Let $\cB$ be an admissible subcategory of a triangulated category $\cA$.
Assume that $\cB ^\perp$ is right admissible.
Then $\cA =(\cB,\cB ^\perp)$ is a semiorthogonal decomposition.
\end{example}

\begin{lemma} Let $\cA=(\cB ,\cC)$ be a semiorthogonal decomposition. Then $\cB=\cC^\perp$.
\end{lemma}

\begin{pf} Let $\cD
\subset \cC^\perp$ be the right orthogonal to $\cB$ in $\cC ^\perp$.  Then $\cD =0$, because
$\cA$ is the triangulated envelope of $\cB \cup \cC$.
By Corollary 1.17 above the category $\cB$ is admissible in $\cC ^\perp$ and hence $\cB=\cC ^\perp$.
\end{pf}

\begin{remark} Let $\cA =(\cB _1,...,\cB _n)$ be a semiorthogonal decomposition. Assume that
for each $i$ there is given a semiorthogonal decomposition $\cB _i=(\cC _i^1,...,\cC _i ^{k_i})$.
Then by Corollary 1.17 each $\cC_i^j$ is admissible in $\cA$ and
$\cA =(\cC _1^1,...,\cC _n^{k_n})$ is a semiorthogonal decomposition.
\end{remark}

\begin{defn} Let $\cT$ be the abelian group generated by equivalence classes of triangulated
categories with the relations coming from semiorthogonal decompositions. Namely,
we put $[\cA]=[\cB]+[\cC]$ if there exists a semiorthogonal decomposition
$\cA _1=(\cB _1,\cB _2)$ with $\cA _1$, $\cB _1$ $\cB _2$ being equivalent to $\cA$,
$\cB$, $\cC$ respectively. We call $\cT$ the {\it Grothendieck group} of triangulated categories.
\end{defn}

This work grew out of an attempt to make $\cT$ into a commutative ring by defining
an appropriate product of triangulated categories.

\section{Tensor product of triangulated categories: a failed attempt}

Recall that there is a simple construction of the tensor product of preadditive (or
prelinear) categories. Namely, if $\cE$ and $\cF$ are two such categories, we let the
objects of $\cE \otimes \cF$ to be the pairs $(X,Y)\in Ob\cE \times Ob\cF$ with
$Mor((X_1,Y_1),(X_2,Y_2))=Mor(X_1,X_2)\otimes Mor(Y_1,Y_2)$.

An axiomatic definition of the tensor product of abelian categories was given by Deligne in
[De]. Let us recall it. Let $\{\cA _i\}_{i\in I}$ be a collection of ($k$-linear)
abelian categories. An abelian category $\cA$ together with a functor
$$\otimes :\prod \cA_i\longrightarrow \cA$$
which is right exact in each variable
is called the {\it tensor product} of $\cA_i$'s if it has the following universal property:
For any abelian category $\cC$ the composition with the functor $\otimes $ establishes an
equivalence of categories of right exact functors $\cA \to \cC$ and functors
$\prod \cA_i\to \cC$, which are right exact in each variable.

Deligne proves the existence of the tensor product under some finiteness condition on
the categories $\cA_i$. Namely, every object should have finite length, and the $\Hom $-space
between any two objects should be finite dimensional. Actually in this case the functor
$\otimes $ is exact in each variable.

One could try to give a similar definition for triangulated categories, say by
replacing right exact functors by exact ones (in the triangulated sense). However this
approach meets difficulties when one tries to prove the
existence, because the triangulated categories are not as rigid as
abelian ones (since taking the cone of a morphism is not a functorial operation).
Also one might want the tensor product of saturated triangulated categories to be
also a saturated triangulated category. In this case a good candidate seems to be
$Fun_{multi-ex}(\cA _1\otimes ...\otimes \cA _n, D^b(Vect))$ -- the category of
multi-exact functors from $\cA _1\otimes ...\otimes \cA _n$ to $D^b(Vect)$. Unfortunately,
this last category does not have a visible triangulated structure.

So we came to conclusion that it is necessary to consider {\it enhanced }
triangulated categories.

\section{DG categories, pretriangulated categories, enhanced triangulated categories}

Our main references here are [BoKa2],[Dr],[Ke].

\subsection{DG categories}
 A DG category is an additive category $\cA$ in which the sets $\Hom (A,B)$, $A,B\in Ob\cA$,
 are
 provided
with a structure of a $\bbZ$-graded $k$-module and a differential $d:\Hom(A,B)\to \Hom (A,B)$
of degree 1, so that for every $A,B,C\in \cA$ the composition $\Hom (A,B)\times \Hom (B,C)
\to \Hom (A,C)$ comes from a morphism of complexes $\Hom (A,B)\otimes \Hom (B,C)
\to \Hom (A,C)$. Also there is a closed degree zero morphism $1_A\in \Hom (A,A)$ which behaves
as the identity under composition of morphisms.

The simplest example of a DG category is the category $DG(k)$ of complexes of
$k$-vector spaces, or DG $k$-modules.

Using the supercommutativity isomorphism $S\otimes T\simeq T\otimes S$ in the
category of DG $k$-modules one defines for every DG category $\cA$ the opposite DG category
$\cA ^0$ with $Ob\cA ^0=Ob\cA$, $\Hom_{\cA ^0}(A,B)=\Hom _{\cA}(B,A)$.
We denote by $\cA ^{gr}$ the {\it graded }
 category which is obtained from $\cA$ by forgetting the differentials on $\Hom $'s.

The tensor product of DG-categories $\cA$ and $\cB$ is defined as follows:

(i) $Ob(\cA \otimes \cB):=Ob\cA \times Ob\cB$; for $A\in Ob\cA$ and $B\in Ob\cB$ the
corresponding object is denoted by $A\otimes B$;

(ii) $\Hom(A\otimes B,A^\prime \otimes B^\prime):=\Hom (A,A^\prime)\otimes \Hom (B,B^\prime)$
and the composition map is defined by $(f_1\otimes g_1)(f_2\otimes g_2):=
(-1)^{\deg(g_1)\deg(f_2)}f_1f_2\otimes g_1g_2.$

Note that the DG categories $\cA \otimes \cB$ and $\cB \otimes \cA$ are canonically
isomorphic. In the above notation the isomorphism functor $\phi$ is
$$\phi (A\otimes B)=(B\otimes A), \quad \phi(f\otimes g)=(-1)^{\deg(f)\deg(g)}(g\otimes f).$$

Given a DG-category $\cA$ one defines the graded category $Ho^\cdot (\cA)$ with
$ObHo^\cdot (\cA)=Ob\cA$ by replacing each $\Hom$ complex by the direct sum of its cohomology
groups. We call $Ho^\cdot (\cA)$ the {\it graded homotopy category} of $\cA$. Restricting
ourselves to the 0-th cohomology of the $\Hom $ complexes we get the {\it homotopy category}
$Ho(\cA)$.

Two objects $A,B\in Ob\cA$ are called DG {\it isomorphic} (or, simply, isomorphic) if there
exists an invertible degree zero morphism $f\in \Hom(A,B)$. We say that $A,B$ are
{\it homotopy equivalent} if they are isomorphic in  $Ho(\cA)$.

A DG-functor between DG-categories $F:\cA \to \cB$ is said to be a
{\it quasi-equivalence} if
$Ho^\cdot(F):Ho^\cdot (\cA)\to Ho^\cdot (\cB)$ is full and faithful and $Ho(F)$ is essentially
surjective.
We say that $F$ is a DG {\it equivalence} if it is full and faithful and
every object of $\cB$ is DG isomorphic to  an object
of $F(\cA)$. Certainly, a DG equivalence is a quasi-equivalence.
DG categories $\cC$ and $\cD$ are called {\it quasi-equivalent} if there exist
DG categories $\cA _1,...,\cA _n$ and a chain of quasi-equivalences
$$\cC \leftarrow \cA _1 \rightarrow ...\leftarrow \cA _n \rightarrow \cD.$$

Given DG categories $\cA$ and $\cB$ the collection of covariant DG functors
$\cA \to \cB$ is itself the collection of objects of a DG category, which we
denote by $Fun _{DG}(\cA ,\cB)$. Namely, let $\phi $ and $\psi$ be two DG functors.
Put $\Hom ^k(\phi ,\psi)$ equal to the set of natural transformations
$t:\phi ^{gr} \to \psi ^{gr}[k]$ of graded functors from $\cA ^{gr}$ to
$\cB ^{gr}$.
This means that for any morphism $f \in \Hom_{\cA}^s(A,B)$ one has
$$\psi (f )\cdot t(A)=(-1)^{ks}t(B)\cdot \phi (f).$$
On each $A\in \cA$ the differential of the transformation $t$ is equal
to $(dt)(A)$ (one easily checks that this is well defined). Thus,
the closed transformations of degree 0 are the
DG transformations of DG functors. A similar definition gives us the DG-category
 consisting of the contravariant DG functors $Fun _{DG}(\cA ^0 ,\cB)=Fun _{DG}(\cA  ,\cB ^0)$
 from $\cA$ to $\cB$.

\begin{remark} Let $\cA$ and $\cB$ be DG categories. Note that the category $Ho(\cA)\otimes
Ho(\cB)$ (and hence $Ho(\cA)\times Ho(\cB)$) is a (not full, in general) subcategory of
$Ho(\cA\otimes \cB)$.
\end{remark}

\subsection{DG modules over DG categories}
We denote the DG category $Fun _{DG}(\cB ,DG(k))$ by $\cB \text{-mod}$ and call it the
 category
of DG $\cB$-modules. There is a natural covariant DG functor
$h:\cA \to \cA ^0\text{-mod}$
(the Yoneda embedding) defined by $h^A(B):=\Hom _{\cA}(B,A)$. As in the "classical" case
one verifies that the functor $h$ is full and faithful, i.e.
$$\Hom _{\cA}(A,A^\prime)=\Hom_{\cA ^0\text{-mod}}(h^A,h^{A^\prime}).$$
Moreover, for any $F\in \cA ^0\text{-mod}$, $A\in \cA$
$$\Hom _{\cA ^0\text{-mod}}(h^A,F)=F(A).$$

The $\cA^0$-DG-modules $h^A$, $A\in \cA$ are called {\it free}. An $\cA ^0$-DG-module
$F$ is called {\it semi-free} if it has a filtration
$$0=F_0\subset F_1\subset ...=F,$$
such that $F_{i+1}/F_i$ is isomorphic to a direct sum of shifted free $\cA^0$-DG-modules
$h^A[n]$, $n\in \bbZ$. The full
subcategory of semi-free $\cA ^0$-DG-modules is denoted by $SF(\cA)$.

An $\cA^0$-DG-module $F$ is called acyclic, if the complex $F(A)$ is acyclic for all
$A\in \cA$. Let $D(\cA)$ denote the {\it derived category} of $\cA ^0$-DG-modules, i.e.
$D(\cA)$ is the Verdier quotient of the homotopy category $Ho(\cA ^0\text{-mod})$ by the
subcategory of acyclic DG-modules. The following proposition was essentially proved in [Ke].

\begin{prop}[Dr]  The inclusion functor $SF(\cA )\hookrightarrow \cA ^0\text{-mod}$
 induces an equivalence of triangulated categories $Ho(SF(\cA))\simeq
 D(\cA)$.
 \end{prop}

A DG functor $G:\cA \to \cB$ induces DG functors
$$\Res _G:\cB ^0\text{-mod}\to \cA ^0\text{-mod}, \quad \quad \Ind_G:\cA ^0\text{-mod}
\to \cB ^0\text{-mod},$$
where $\Res _G$ is the obvious restriction functor. We refer to [Dr] for the definition of
$\Ind _G$; let us only list some of the properties.

\medskip

\noindent1). The functor  $\Ind _G$ is left adjoint to $\Res _G$, that is for every $\Phi \in
\cA ^0\text{-mod}$, $\Psi \in \cB ^0\text{-mod}$ there is a canonical isomorphism of complexes
$$\Hom _{\cB ^0\text{-mod}}(\Ind _G(\Phi),\Psi )=\Hom _{ \cA ^0\text{-mod}}(\Phi ,\Res
_G(\Psi)).$$

\noindent2). For any $A\in \cA$, $\Ind _G(h^A)=h^{(G(A))}$, which means that the following
diagram is commutative
$$\begin{array}{ccc}
\cA & \stackrel{h}{\longrightarrow} & \cA ^0\text{-mod}\\
G \downarrow & & \downarrow \Ind _G\\
\cB & \stackrel{h}{\longrightarrow} & \cB ^0\text{-mod},
\end{array}
$$
where the horizontal arrows are the Yoneda embeddings.

\noindent3). The functor $\Ind _G$ preserves semi-free DG modules and $\Ind _G:SF(\cA)\to
SF(\cB)$ is a quasi-equivalence if $G$ is such.

\medskip

\begin{cor} [Ke] If DG categories $\cA$ and $\cB$ are quasi-equivalent, then the
derived categories $D(\cA)$ and $D(\cB)$ are equivalent.
\end{cor}

\begin{pf} Indeed, this follows from Proposition 3.2 and the last property of the
functor $\Ind _G$.
\end{pf}

It is natural to consider the category of $h$-projective $\cA ^0$-DG-modules
which we introduce next.  Namely, we call a  $\cA ^0$-DG-module $P$ $h$-{\it projective} if
$$\Hom_{Ho(\cA ^0\text{-mod})}(P,F)=0$$
for every acyclic $F\in \cA^0\text{-mod}$. Let  $\cP(\cA)\subset \cA ^0\text{-mod}$ denote the full
subcategory of $h$-projective objects. It can be shown that a semi-free
$\cA ^0$-DG-module is $h$-projective. Hence the Proposition 3.2 implies the equivalences
$$Ho(SF(\cA))\simeq Ho(\cP(\cA))\simeq D(\cA).$$
Property 3 of the functor $\Ind _G$ implies the following corollary.

\begin{cor} If DG categories $\cA$ and $\cB$ are quasi-equivalent, then also $\cP(\cA)$ and
$\cP(\cB)$ are quasi-equivalent.
\end{cor}

\subsection{Pretriangulated DG categories}
   Given a DG-category $\cA$ one can associate to it a triangulated category $\cA ^{tr}$
   [BoKa2]. It is defined as the homotopy category of a certain DG-category $\cA^{\text{pre-tr}}$.
   The idea of the definition of $\cA ^{\text{pre-tr}}$ is to formally add cones of all morphisms,
   cones of morphisms between cones, etc.

   First we need to clarify the notion of a "formal shift" of an object.

   \begin{defn} We define the DG category $\bar {\cA}$ as follows:
   $$Ob\bar{\cA}=\{A[n]\vert A\in Ob\cA, \ n\in \bbZ\},$$
   and
   $$\Hom _{\bar{\cA}}(A[k],B[n])=\Hom _{\cA}(A,B)[n-k]$$
   as {\it graded vector spaces}. If $f\in \Hom _{\cA}(A,B)$ is considered as an element of
 $\Hom _{\bar{\cA}}(A[k],B[n])$
   under the above identification  then the  differentials are related by the formula
   $$d_{\bar{\cA}}(f)=(-1)^nd_{\cA}(f).$$
\end{defn}

Notice, for example, that the differential in $\Hom _{\bar{\cA}}(A[1],B[1])$ is
equal to {\it minus} the differential in $\Hom _{\cA}(A,B)$.

One can check that the composition of morphisms in  $\bar {\cA}$ is compatible with
the Leibniz rule, so that $\bar {\cA}$ is indeed a DG category.
Clearly  $\bar{\cA}$ contains $\cA$ as a full DG subcategory.

Given an object $A\in \cA$ the object $A[r]$ is characterized
(up to a DG isomorphism)
by the existence of closed morphisms $f:A\to A[r]$, $g:A[r]\to A$ of degrees $-r$ and $r$
respectively, such that $fg=gf=1$. Thus in particular every DG functor commutes with shifts.

 \begin{defn}  The objects of $\cA^{\text{pre-tr}}$ are "one-sided
    twisted complexes", i.e. formal expressions $(\oplus_{i=1}^nC_i[r_i],q)$, where
    $C_i\in Ob\cA$, $r_i\in \bbZ$, $n\geq 0$, $q=(q_{ij})$, $q_{ij}\in
    \Hom(C_j[r_j],C_i[r_i])$
    is homogeneous of degree 1, $q_{ij}=0$ for $i\geq j,$ $dq+q^2=0$. If $C,C^\prime \in Ob
    \cA ^{\text{pre-tr}}$, $C=(\oplus C_j[r_j],q),$ $C^\prime =(\oplus
    C^\prime_j[r^\prime _j],q^\prime),$ the  $\bbZ$-graded $k$-module $\Hom (C,C^\prime)$
    is the space of matrices $f=(f_{ij})$, $f_{ij}\in \Hom (C_j[r_j],C^\prime _i[r^\prime_i])$,
    and the composition map $\Hom(C,C^\prime)\otimes
    \Hom (C^\prime,C^{\prime \prime})\to \Hom (C,C^{\prime \prime})$ is matrix multiplication.
    The differential $d:\Hom (C,C^\prime)\to \Hom (C,C^\prime )$ is defined by $df:=(df_{ij})+
    q^\prime f-(-1)^lfq$ if $\deg f_{ij}=l$.
    \end{defn}

Notice that  the DG category $\cA^{\text{pre-tr}}$
is closed under formal shifts:
$$(\oplus_{i=1}^nC_i[r_i],q)[1]=(\oplus_{i=1}^nC_i[r_i+1], -q).$$

\begin{defn} Let $\cB$ be a DG category and $f\in \Hom (A,B)$ be a closed degree zero
morphism in
$\cB$. An object $C\in \cB$ is called the cone of $f$, denoted $Cone(f)$, if $\cB$ contains
the object $A[1]$ and there exist degree zero morphisms
$$A[1]\stackrel{i}{\rightarrow} C\stackrel{p}{\rightarrow} A[1], \quad
B\stackrel{j}{\rightarrow} C\stackrel{s}{\rightarrow} B$$
with the properties
$$p i=1, \ \ s j=1,\ \ s i=0,\ \ p j=0, \ \  i p+j s=1,$$
and
$$d(j)=d(p)=0,\ \ d(i)=j f,\ \ d(s)=-f p.$$
\end{defn}

\begin{lemma} The cone of a closed degree zero morphism is uniquely defined up to a
DG isomorphism.
\end{lemma}

\begin{pf} Note that the first set of conditions means that $C$ is the direct sum of $A[1]$ and $B$
in the corresponding graded category $\cB ^{gr}$. Thus for any object $E$ in $\cA$ there are
isomorphisms of graded $k$-modules
$$\Hom (E,C)=\Hom (E,A[1])\oplus \Hom (E,B),$$
$$ \Hom (C,E)=\Hom (A[1],E)\oplus \Hom (B,E),$$
which are given by composing with $i$ and $j$ (or with $p$ and $s$). Then the second
set of conditions determines the differentials in $\Hom (E,C)$ and $\Hom (C,E)$.
\end{pf}

Given a closed degree zero morphism $f:A\to B$ the diagram
$$A\stackrel {f}{\longrightarrow} B\stackrel{j}{\longrightarrow} Cone(f)
\stackrel{p}{\longrightarrow} A[1]$$
is called a preexact triangle.

\begin{remark} It is clear that any DG functor preserves cones of closed degree zero morphisms
and preserves preexact triangles.
\end{remark}

\begin{prop}[BoKa2] Let $\cA$ be a DG category. Then

a) the DG category $\cA ^{\text{pre-tr}}$ is closed under taking cones of closed degree zero
morphisms;

b) every object in $\cA ^{\text{pre-tr}}$ can be obtained from objects in $\cA$ by
taking successive cones of closed degree zero morphisms.
\end{prop}

\begin{pf} a). Given a closed morphism of degree zero
$$f:(\oplus C_i[r_i],q)\to (\oplus
    C^\prime_j[r^\prime _j],q^\prime)$$
its cone is the twisted complex $(\oplus
 C^\prime_j[r^\prime _j]\oplus C_i[r_i+1],(q^\prime,-q+f))$. For example,
if $A,B\in \cA$ and $f:A\to B$
    is a closed morphism of degree 0 then $Cone(f)$ is the twisted complex
    $(B\oplus A[1] ,(0,f))\in
    \cA ^{\text{pre-tr}}$.

b). Let $C=(\oplus_{i=1}^nC_i[r_i],q)$ be a twisted complex.
Consider its twisted subcomplex $C^\prime =(\oplus_{i=1}^{n-1}C_i[r_i],q^\prime)$,
where $q^\prime =q-\oplus _{i}q_{in}$.
Then $C$ is the cone of
the closed degree zero morphism $\oplus_{i=1}^{n-1}q_{in}:(C_n[r_n-1],0)\to C^\prime$.
\end{pf}

A DG category $\cA$ is said to be {\it pretriangulated} if for every
$A\in \cA$, $k\in \bbZ$ the object
$A[k]\in \cA ^{\text{pre-tr}}$ is homotopy equivalent to an object of $\cA$ and for every closed
morphism of degree zero $f$ in $\cA$ the object $Cone(f)\in \cA ^{\text{pre-tr}}$ is homotopy
equivalent to an object of $\cA$. We say that $\cA$ is {\it strongly pretriangulated} if the
same is true with "homotopy equivalent" replaced by "DG isomorphic". Actually, if $\cA$
is pretriangulated (resp. strongly pretriangulated) then every object of $\cA ^{\text{pre-tr}}$ is
homotopy equivalent (resp. DG isomorphic) to an object of $\cA$ [Dr]. Thus $\cA$ is
pretriangulated (resp. strongly pretriangulated) if and only if the embedding
 $Ho(\cA)\hookrightarrow Ho(\cA^{\text{pre-tr}})$ is an equivalence (resp. the embedding
$\cA \hookrightarrow \cA ^{\text{pre-tr}}$ is a DG equivalence).

\begin{prop} Let $\cA$, $\cB$ be DG categories and $F:\cA \to \cB$ be a DG functor.
Assume that $\cB$ is strongly pretriangulated. Then there
exists a unique (up to a DG isomorphism)
DG functor $G:\cA ^{\text{pre-tr}}\to \cB$ which extends $F$.
\end{prop}

\begin{pf} The canonical embedding $J: \cB \to
\cB ^{\text{pre-tr}}$ is a DG equivalence. Thus we may put $G=J^{-1}\cdot
F ^{\text{pre-tr}}$. The uniqueness of $G$ follows again from part b) of Proposition 3.10.
\end{pf}

For a DG-category $\cA$ the DG categories $\cA ^{\text{pre-tr}}$, $\cA ^0\text{-mod}$, $SF(\cA)$,
$\cP(\cA)$ are strongly pretriangulated.
Moreover, if $\cA$ is a DG category and $\cB$ is a strongly pretriangulated category, then
the DG category $Fun_{DG}(\cA ,\cB)$ is strongly pretriangulated.

    If $\cA$ is \text{pretr}iangulated then every closed degree 0 morphism $f:A\to B$ in
    $\cA$ gives rise to the usual triangle $A\to B\to Cone(f)\to A[1]$ in $Ho (\cA)$.
    Triangles of this type and those isomorphic to them are called {\it exact}. If $\cA$
    is pretriangulated then $Ho (\cA)$ becomes a triangulated category. So
    $\cA^{tr}:=Ho(\cA ^{\text{pre-tr}})$ is a triangulated category.

Let $F:\cA \hookrightarrow \cA ^{\text{pre-tr}}$ denote (temporarily)
the canonical (full and faithful)
embedding. We have the commutative diagrams
$$\begin{array}{ccc}
\cA & \stackrel{h_{\cA}}{\longrightarrow} & \cA ^0\text{-mod}\\
F\downarrow & & \downarrow \Ind _F \\
\cA ^{\text{pre-tr}} & \stackrel{h_{\cA ^{\text{pre-tr}}}}{\longrightarrow} &
(\cA ^{\text{pre-tr}})^0\text{-mod}
\end{array}
$$
and
$$\begin{array}{ccc}
\cA & \stackrel{h_{\cA}}{\longrightarrow} & \cA ^0\text{-mod}\\
F\downarrow & & \uparrow \Res _F \\
\cA ^{\text{pre-tr}} & \stackrel{h_{\cA ^{\text{pre-tr}}}}{\longrightarrow} &
(\cA ^{\text{pre-tr}})^0\text{-mod}
\end{array}
$$
The DG functors $\Ind _F$ and $\Res _F$ are mutually inverse DG
equivalences, which preserve the subcategories of semi-free DG modules. We denote
the full and faithful DG functor
$$\alpha =\alpha _{\cA}:=\Res _F\cdot h_{\cA ^{\text{pre-tr}}}:\cA ^{\text{pre-tr}}
\to \cA ^0\text{-mod}.$$
Then $\alpha (\cA ^{\text{pre-tr}})$ is a strictly full subcategory of $SF(\cA)$.

\begin{remark}[Dr]
    If a DG functor $G:\cA \to \cB$ is a quasi-equivalence, then the
corresponding quasi-equivalence $\Ind _G:SF(\cA)\to SF(\cB)$ induces a quasi-equivalence
$$\Ind _G:\alpha (\cA ^{\text{pre-tr}})\to \alpha (\cB ^{\text{pre-tr}}).$$
\end{remark}

\subsection{Perfect DG modules}

\begin{defn}  Let $\cA$ be a DG category.
Consider the full pretriangulated subcategory $\alpha (\cA ^{\text{pre-tr}})\subset SF(\cA)$.
Let $\text{Perf-}\cA$ be the full DG subcategory of $SF(\cA)$ consisting of DG modules which
are homotopy equivalent to a direct summand of an object in $\alpha (\cA ^{\text{pre-tr}})$.
We call these the {\it perfect} $\cA^0$-DG-modules.
\end{defn}

The full subcategory $Ho(\text{Perf-}\cA)\subset Ho(SF(\cA))$ is the epaisse envelope of the
triangulated subcategory $Ho(\alpha (\cA ^{\text{pre-tr}}))\subset Ho(SF(\cA))$. Hence
$Ho(\text{Perf-}\cA)$ is also triangulated and, therefore, $\text{Perf-}\cA$ is a
(strongly) pretriangulated
subcategory of $SF(\cA)$.
Note that an arbitrary
direct sum of semi-free $\cA ^0$-DG-modules is again semi-free. Hence the
triangulated category $Ho(SF(\cA))$ contains arbitrary direct sums. By Theorem 1.11 above
it is Karoubian. Hence the category $Ho(\text{Perf-}\cA)$ is also Karoubian, and thus it is
the Karoubization of the triangulated category $\cA ^{tr}$.

It was pointed to us by Bernhard Keller that the category $Ho(\text{Perf-}\cA)$
 can be characterized as consisting of {\it compact} or {\it small} objects in
 $Ho(SF(\cA))\simeq D(\cA)$ [Ke],[Ne],[Ra].

Note that the Yoneda embedding $h:\cA \to \cA ^0\text{-mod}$ defines a full and faithful
DG functor $h:\cA \to \text{Perf-}\cA$. Also a DG functor $F:\cA \to \cB$ induces a DG functor
$\Ind _F:\text{Perf-}\cA \to \text{Perf-}\cB$ so that the functorial diagram
$$\begin{array}{ccc}
\cA & \stackrel{h}{\longrightarrow} & \text{Perf-}\cA\\
F \downarrow & & \downarrow \Ind_F \\
\cB & \stackrel{h}{\longrightarrow} & \text{Perf-}\cB
\end{array}
$$
is commutative.

\begin{lemma} If $F:\cA\to \cB$ is a quasi-equivalence, then $\Ind _F: \text{Perf-}\cA
\to \text{Perf-}\cB$ is also such.
\end{lemma}

\begin{pf} The quasi-equivalence $F:\cA \to \cB$ induces a quasi-equivalence
 $\Ind _F:\alpha _{\cA}(\cA ^{\text{pre-tr}})
\to \alpha _{\cB}(\cB ^{\text{pre-tr}})$. Since the subcategories $Ho(\text{Perf-}\cA)\subset
Ho(SF(\cA))$ and
$Ho(\text{Perf-}\cB)\subset Ho(SF(\cB))$ are the epaisse envelopes of the subcategories
$Ho(\alpha _{\cA}(\cA ^{\text{pre-tr}}))$ and $Ho(\alpha _{\cB}(\cB ^{\text{pre-tr}}))$ respectively, the
DG functor $\Ind_F$ induces the quasi-equivalence $\text{Perf-}\cA \to \text{Perf-}\cB$.
\end{pf}

\begin{cor} If DG categories $\cA $ and $\cB$ are quasi-equivalent, then so are
$\text{Perf-}\cA $ and $\text{Perf-}\cB$.
\end{cor}

\begin{pf} This follows from the last lemma.
\end{pf}

\begin{lemma} Let $\cA$ be a DG category, $\cB\subset \cA$ a full DG subcategory such that
$Ho(\cB)$ is dense in $Ho(\cA)$, i.e. every object $A\in Ho(\cA)$ is a direct summand of an
object in $Ho(\cB)$. Then $Ho(\cB ^{\text{pre-tr}})$ is dense in $Ho(\cA ^{\text{pre-tr}})$.
\end{lemma}

\begin{pf} Let $A,B,C,D\in \cA$ and assume that $A\oplus B$ and $C\oplus D$ are homotopy
equivalent to objects of $\cB$. Let $f:A\to C$ be a closed morphism of degree zero.
By part b) of Proposition 3.9 it suffices
to show that there exists $K\in \cA^{\text{pre-tr}}$ such that
$Cone(f)\oplus K$ is homotopy equivalent to
an object of $\cB ^{\text{pre-tr}}$. Consider the closed morphism of degree zero
$g:A\oplus B\to C\oplus D$, where $g\vert _A=f$ and $g\vert _B=0$. Then
$$Cone(g)=Cone(f)\oplus D\oplus B[1],$$
which is homotopy equivalent to an object of $\cB ^{\text{pre-tr}}$, since $A\oplus B$,
and $C\oplus D$ are such.
\end{pf}

\begin{prop} Let $\cA$ and $\cB$ be DG categories. The natural full and faithful DG functor
$G:\cA \otimes \cB\to (\text{Perf-}\cA)\otimes \cB$ induces a quasi-equivalence
$\Ind _G:\text{Perf-}(\cA \otimes \cB)\to \text{Perf-}((\text{Perf-}\cA)\otimes \cB)$.
\end{prop}

\begin{pf} Denote $\cC=\cA \otimes \cB$, $\cD = (\text{Perf-}\cA)\otimes \cB$. Consider the
induced full and faithful DG functor
$$\Ind _G:\alpha _{\cC}(\cC ^{\text{pre-tr}})\to \alpha _{\cD}(\cD ^{\text{pre-tr}}).$$
Since $\cC ^{\text{pre-tr}}$ is DG equivalent to
$(\cA ^{\text{pre-tr}}\otimes \cB)^{\text{pre-tr}}$  and $Ho(\cA ^{\text{pre-tr}}\otimes \cB)$
is dense in $Ho((\text{Perf-}\cA)\otimes \cB)$ it follows from Lemma 3.15 that
 the category $Ho(\Ind _G)(Ho(\alpha _{\cC}(\cC ^{\text{pre-tr}})))$ is dense in
$Ho(\alpha _{\cD}(\cD ^{\text{pre-tr}}))$. Therefore $Ho(\Ind _G)$ induces an equivalence
of Karoubizations $Ho(\text{Perf-}\cC)$ and  $Ho(\text{Perf-}\cD)$ of these categories.
\end{pf}

\begin{defn} A DG category $\cA$ is {\it perfect} if it is pretriangulated and the
triangulated category $Ho(\cA)$ is Karoubian.
\end{defn}

Note that a DG category which is quasi-equivalent to a perfect one is by itself perfect.

\begin{example} $\text{Perf-}\cA$ is perfect for any DG category $\cA$.
\end{example}

\begin{prop} If a DG category $\cA$ is perfect, then the (Yoneda)
 embedding $\cA \hookrightarrow
\text{Perf-}\cA$ is a quasi-equivalence. In particular $\text{Perf-}\cA$ is quasi-equivalent
to $\text{Perf-}(\text{Perf-}\cA)$.
\end{prop}

\begin{pf} Similar to the proof of Proposition 3.16
\end{pf}

\subsection{Enhanced triangulated categories}
    Given a triangulated category $\cD$, by its {\it enhancement} we shall mean a
    \text{pre-tr}iangulated DG category $\cA$ together with an equivalence of triangulated
    categories $\epsilon _{\cA}:Ho(\cA)\to \cD$.
    The category $\cD$ is then said to be {\it enhanced}.
Given enhancements $(\cA ,\epsilon _{\cA})$ and $(\cB ,\epsilon _{\cB})$ of $\cD$ we call
a DG functor $F:\cA \to \cB$ a {\it quasi-equivalence of enhancements} if $F$ is a
quasi-equivalence and the functors $\epsilon _{\cA}$ and $Ho(F)\circ \epsilon _{\cB}$ are
isomorphic. Enhancements $(\cA ,\epsilon _{\cA})$ and $(\cB ,\epsilon _{\cB})$ are called
quasi-equivalent if there exist enhancements $(\cA _i,\epsilon _{\cA _i})$ of $\cD$ and
 a chain of quasi-equivalences of enhancements
$$\cA \leftarrow \cA _1 \rightarrow \cA _2 ...\leftarrow \cA _n \rightarrow \cB.$$
It is natural to consider quasi-equivalence classes of enhancements of a
triangulated category.

\section{Grothendieck ring of pretriangulated categories}

\subsection{Grothendieck group of pretriangulated categories}

\begin{defn} Denote by $\cP\cT$ the abelian group generated by quasi-equivalence classes
of pretriangulated categories with relations coming from semi-orthogonal decompositions.
Namely, given pretriangulated categories $\cA$, $\cB$, $\cC$, there is the relation
$$[\cA]=[\cB]+[\cC]$$
in $\cP\cT$ if there exist pretriangulated categories $\cA ^\prime$, $\cB ^\prime$,
$\cC ^\prime$ which are quasi-equivalent to $\cA$, $\cB$, $\cC$ respectively, such that

1) $\cB ^\prime$, $\cC ^\prime$ are  DG subcategories of $\cA ^\prime$,

2) $Ho(\cB ^\prime)$, $Ho(\cC ^\prime)$ are admissible subcategories of $Ho(\cA ^\prime)$, and

3) $Ho(\cA ^\prime)=(Ho(\cB ^\prime), Ho(\cC ^\prime))$ is a semiorthogonal decomposition.
\end{defn}

We are going to turn the group $\cP\cT$ into an associative commutative ring by defining
the appropriate product of pretriangulated categories.

\subsection{$\bullet$ - product of pretriangulated categories}

\begin{defn} Let $\cA _1$, ..., $\cA _n$ be DG categories. Define the DG category
$$\cA _1\bullet ...\bullet \cA _n:= \text{Perf-}(\cA _1\otimes ...\otimes \cA _n).$$
Thus $\cA _1\bullet ...\bullet \cA _n$ is a perfect DG category and hence the category
$Ho(\cA _1\bullet ...\bullet \cA _n)$ is triangulated and Karoubian.
\end{defn}

\begin{remark} The product $\cA _1\bullet ...\bullet \cA _n$ is functorial in each variable
with respect to DG functors.
\end{remark}

Proposition 3.14 above implies that $\bullet $ preserves quasi-equivalence classes of
DG categories. Proposition 3.16 implies that the induced operation
$$\bullet :\cP\cT \times \cP \cT \to \cP \cT$$
is associative. Since the DG categories $\cA\otimes \cB$ and $\cB\otimes \cA$ are isomorphic,
this operation is commutative. We want to show that it is distributive with respect to the
addition in the group $\cP\cT$.

 \begin{lemma} If a DG functor $F:\cA _1 \to \cA _2$ is such that $Ho(F)$ is full and faithful,
 then
 the functor $Ho(\Ind _F):Ho(\text{Perf-}\cA _1)\to Ho(\text{Perf-}\cA _2)$ is also full
 and faithful.
 \end{lemma}

 \begin{pf} Indeed, it follows from part b) of Proposition 3.9
  that the exact functor $Ho(\Ind_F):\cA _1^{tr}\to \cA _2^{tr}$
 is full and faithful, hence so is the functor
 $Ho(\Ind _F):Ho(\text{Perf-}\cA _1)\to Ho(\text{Perf-}\cA _2)$.
 \end{pf}

\begin{lemma} Let $\cA$, $\cD$ be  pretriangulated DG categories and
$\cB \subset \cA$ be a pretriangulated DG subcategory such that
$Ho(\cB)$  is admissible in $Ho(\cA)$. Then  $\cB \bullet \cD$ is
pretriangulated DG subcategory of  $\cA \bullet \cD$ such that $Ho(\cB \bullet
\cD)$ is admissible in $Ho( \cA \bullet \cD)$.
\end{lemma}

\begin{pf} Assume for example that $Ho(\cB)$ is right admissible in
$Ho(\cA)$. There exists a pretriangulated DG subcategory
$\cC\subset \cA$ such that $Ho(\cC)=Ho(\cB)^\perp$. Fix one such subcategory $\cC$.

 Denote by $i:\cB \otimes \cD \hookrightarrow \cA \otimes \cD$ and
$j:\cC \otimes \cD \hookrightarrow \cA \otimes \cD$ the embedding
DG functors.
The functors $Ho(i)$ and $Ho(j)$ are full and faithful, so
by the last lemma we may identify $Ho(\cB \bullet \cD)$ and $Ho(\cC \bullet \cD)$
as full subcategories of $Ho(\cA \bullet \cD)$, which we denote $\tilde{\cB}$ and
$\tilde{\cC}$ respectively.

Note that $\tilde{\cC}\subset \tilde{\cB}^\perp$. It suffices to prove that
for every $E\in Ho(\cA \bullet \cD)$ there exists an exact triangle
$$E_{\tilde{\cB}}\to E\to E_{\tilde{\cC}},$$
where $E_{\tilde{\cB}}\in \tilde{\cB}$ and $E_{\tilde{\cC}}\in \tilde{\cC}$.

Notice that the collection of objects of the subcategory
$Ho(\cA \otimes \cD)\subset Ho(\cA \bullet \cD)$
classically generates the category $Ho(\cA \bullet \cD)$. On the other hand for each object
$A\otimes D\in Ho(\cA \otimes \cD)$ there  exists an
exact triangle
$$A_{\cC}\otimes D\to A\otimes D\to A_{\cB}\otimes D$$
with $A_{\cB}\in Ho(\cB)$ and $A_{\cC}\in Ho(\cC)$. It remains to apply the Lemma 1.20.
The case of $Ho(\cB)$ being left admissible is treated similarly.
\end{pf}

The proof of the last lemma also contains a proof of the next proposition.

 \begin{prop} Let $\cA$, $\cD$ be pretriangulated DG categories and
$\cB,\cC \subset \cA$ be two  DG subcategories
such that $Ho(\cB)$ and $Ho(\cC)$ are admissible in $Ho(\cA)$ and
$Ho(\cA)=(Ho(\cB),Ho(\cC))$ is a semi-orthogonal decomposition. Then $\cB \bullet \cD$,
$\cC\bullet \cD$ are DG subcategories in $\cA \bullet \cD$, such that
$Ho(\cB \bullet \cD)$ and $Ho(\cC\bullet \cD)$ are admissible in $Ho(\cA \bullet \cD)$ and
 $Ho(\cA \bullet \cD)=(Ho(\cB \bullet \cD), Ho(\cC\bullet \cD))$ is a semi-orthogonal
 decomposition.
 \end{prop}

We get the immediate corollary

\begin{cor} The operation $\bullet$ is distributive with respect to addition in the
group $\cP\cT$. So $\cP\cT$ is a commutative associative ring.
\end{cor}

Let $DG(k)^f$ denote the DG category of finite dimensional complexes of $k$-vector spaces.
 Any pretriangulated DG category $\cA$ is quasi-equivalent to $\cA \otimes DG(k)^f$.

\begin{remark} Denote by $\cP\cT ^+\subset \cP\cT$ the subgroup generated by quasi-equivalence
 classes of perfect DG categories. (Note that if $\cA$ is perfect and $\cB \subset \cA$ is a
 pretriangulated subcategory such that $Ho(\cB)$ is admissible in $Ho(\cA)$, then
 $\cB$ is also perfect.) It is clear that $\cP\cT ^+$ is a subring of $\cP\cT$. Actually
$\cP\cT ^+$ is a unital ring with the unit $[DG(k)^f]$.
\end{remark}

\begin{remark} Actually one can make $\cP\cT$ into a ring by using a simpler multiplication
$\cA _1\circ ...\circ\cA _n=(\cA _1\otimes ...\otimes\cA _n)^{\text{pre-tr}}.$ We chose the operation
$\bullet$ because we like triangulated categories which are Karoubian (it gives them a chance
to be saturated).
\end{remark}

    \section{A geometric example}

    \subsection{Standard enhancements of $D(X)$}
 Let $X$ be a smooth projective variety. Let us consider the following model for
 the triangulated category $D(X)$. Consider the abelian category $\cM od (\cO _X)$ of all
 (not necessarily
 quasi-coherent) $\cO _X$-modules. Put $D(X)=D^b_{coh_X}(\cM od _{\cO _X})$. So the objects
 of $D(X)$ and complexes of $\cO _X$-modules with bounded coherent cohomology.

 There are several natural  enhancements of the
 category $D(X)$. Consider the pretriangulated DG-category
 $C(X)$ consisting of bounded below complexes of $\cO _X$-modules with bounded
 coherent cohomology. Let
 $I(X)$ denote the full pretriangulated subcategory of $C(X)$ consisting of complexes
 of injective $\cO _X$-modules. Denote the composition of the natural functors
 $Ho(I(X))\to Ho(C(X))\to D(X)$ by $\epsilon _{I(X)}$.
  It is well known that $\epsilon _{I(X)}$ is an equivalence.

\begin{defn} Let $(\cA ,\epsilon _{\cA})$ be an enhancement of $D(X)$ which is quasi-equivalent
to $(I(X),\epsilon _{I(X)})$. We call $(\cA ,\epsilon _{\cA})$ a {\it standard} enhancement of
$D(X)$.
\end{defn}

  In fact we believe that all enhancements of
 $D(X)$ are standard.

The next lemma is in this direction.

\begin{lemma} Let $\cA\subset C(X)$ be a full pretriangulated DG subcategory such that the
natural functor
$\epsilon _{\cA}:Ho(\cA)\to D(X)$ is an equivalence. Then $(\cA,\epsilon _{\cA})$ is
a standard enhancement of $D(X)$.
\end{lemma}

\begin{pf} Consider the following DG category $\cC$. Objects of $\cC$ are triples
$(A,I,\gamma)$, where $A\in \cA$, $I\in I(X)$, and $\gamma :A\to I$ is an injective
quasi-isomorphism of complexes. The complex $\Hom _{\cC}((A,I,\gamma ),(B,J,\delta ))$
is a subcomplex of $\Hom _{I(X)}(I,J)$ consisting of morphisms which map $A$ to $B$.
We have two obvious DG functors
$$\cA \stackrel{\phi}{\longleftarrow}\cC \stackrel{\psi}{\longrightarrow}I(X),\quad
A\gets (A,I,\gamma)\to I.$$

\noindent{\it Claim}. The DG functors $\phi$ and $\psi$ are quasi-equivalences.

The functor $\phi$ is surjective on objects since every complex in $C(X)$ has a
(bounded below) injective resolution. The functor $Ho(\psi)$ is essentially surjecive,
because $Ho(\cA)$ is equivalent to $D(X)$. So it remains to prove that the functors
$Ho(\phi)$ and $Ho(\psi)$ are full and  faithful. Fix $S=(A,I,\gamma)$, $T=(B,J,
\delta)$ and consider the commutative diagram
$$\begin{array}{rcl}
\Hom _{\cC}(S,T) & \stackrel{\phi _{S,T}}{\longrightarrow} & \Hom_{\cA}(A,B)\\
\psi _{S,T}\downarrow & &\downarrow \delta _*\\
\Hom _{I(X)}(I,J) & \stackrel{\gamma ^* }{\longrightarrow} & \Hom _{C(X)}(A,J).
\end{array}
$$
The map $\gamma ^*$ is a quasi-isomorphism, since $J$ is injective. The map $\delta _*$
is a quasi-isomorphism, because $Ho(\cA)$ is equivalent to $D(X)$. So it suffices to
prove that $\phi _{S,T}$ is a quasi-isomorphism.

Since the complex $J$ is injective the map $\phi _{S,T}$ is surjective and surjective on
cycles. Let $f\in \Hom _{\cC}(S,T)$ be a cycle such that $\phi _{S,T}f=dg$ for some
$g\in \Hom _{\cA}(A,B)$. It suffices to prove that $f$ is a boundary in $\Hom _{\cC}(S,T)$.
Choose $\tilde{g}\in \Hom _{\cC}(S,T)$ such that $\phi _{S,T}\tilde{g}=g$. Replacing
$f$ by $f-d\tilde{g}$ we may assume that $\phi _{S,T}f=0$. But then $f$ is a cycle in
$\Hom_{C(X)}(I/A,J)$. Since $I/A$ is acyclic and $J$ is injective, $f$ is a boundary in
$\Hom _{C(X)}(I/A,J)$; hence also in $\Hom _{\cC}(S,T)$.
 This proves the claim.

 To prove the lemma it suffices to show that the functors $\epsilon _{\cA}\cdot Ho(\phi)$ and
 $\epsilon _{I(X)}\cdot Ho(\psi)$ from $Ho(\cC)$ to $D(X)$ are isomorphic. Fix an object
$(A,I,\gamma)\in \cC$. Then $\epsilon _{\cA}\cdot Ho(\phi)((A,I,\gamma))=A$,
  $\epsilon _{I(X)}\cdot Ho(\psi)((A,I,\gamma))=I$ and the required quasi-isomorphism
  is the embedding $\gamma :A\to I$.
 \end{pf}

\begin{defn} Fix a full subcategory $\cC \subset D(X)$. Denote by $I(\cC)\subset I(X)$
the full DG
subcategory consisting of complexes, which are quasi-isomorphic to objects in $\cC$.
Clearly, the obvious functor $\epsilon _{I(\cC)}:Ho(I(\cC))\to \cC$ is an equivalence.
An enhancement $(\cA,\epsilon _{\cA})$ of $\cC$ is called {\it standard} if it is
quasi-equivalent to $(I(\cC),\epsilon _{I(\cC)})$. Notice that this notion depends
not only on the category $\cC$, but also on the given embedding $\cC \subset D(X)$.
Sometimes we will refer to $cA$ alone as a standard enhancement of $\cC$, if it is clear
what the functor $\epsilon _{\cA}$ is.
\end{defn}

\begin{lemma} Let $\cB \subset C(X)$ be a full DG subcategory such that the
corresponding functor $\epsilon _{\cB}:Ho(\cB)\to D(X)$ is full and faithful.
Denote by $\cC \subset D(X)$ the image of the functor $\epsilon _{\cB}$. Then
$(\cB ,\epsilon _{\cB})$ is a standard enhancement of $\cC$.
\end{lemma}

\begin{pf} Same as that of last lemma.
\end{pf}

The main result of this section is the following theorem.

\begin{thm} Let $X_1, ...,X_n$ be smooth projective varieties over $k$. Then the DG category
$I(X_1)\bullet ...\bullet I(X_n)$ is quasi-isomorphic to $I(X_1\times ...\times X_n)$, i.e.
$$[I(X_1)]\bullet ...\bullet [I(X_n)]=[I(X_1\times ...\times X_n)].$$
\end{thm}

Since the operation $\bullet $ is associative up to quasi-equivalence it suffices to
prove the theorem for $n=2$. Put $X_1=X$, $X_2=Y$.

For a proof of the theorem it will be convenient to use a standard enhancement of $D(X)$
which we presently
define. Let $\cP(X)$ denote the full subcategory of $D(X)$ consisting of perfect complexes,
i.e. finite complexes of vector bundles.
It is well known that every object in $D(X)$ is isomorphic to an object in $\cP(X)$. Hence
the embedding $\cP(X)\subset D(X)$ is an equivalence. Choose a finite affine covering
$\cU$ of $X$. For any $P\in \cP(X)$ consider its $\check C$ech resolution
$P\to C_{\cU}(P)$ defined by $\cU$. Thus $C_{\cU}(P)$ is a finite complex of quasi-coherent
sheaves which are direct sums of sheaves $P^j_U:=i_*i^*P^j$, where $i:U\hookrightarrow X$ is the
open embedding of an affine open subset $U$, which is the intersection of some elements from
$\cU$. Let $\cP(\cU)\subset C(X)$ denote the minimal full DG subcategory
which contains all complexes $C_{\cU}(P)$ for
$P\in \cP(X)$ and is closed under taking cones of closed morphisms of degree zero. Thus
$\cP(\cU)$ is (strongly) pretriangulated.
 (We could denote the DG category $\cP(\cU)$ by  $\cP^+(\cU)$ (resp. the complex
$C_{\cU}(P)$ by $C_{\cU}^+(P)$) as later (proof of Lemma 6.4)
we will also consider the "dual" enhancement
$\cP ^-(\cU)$ using the left $\check C$ech resolutions $C^-_{\cU}(P)\to P$.)

We have the
obvious functor
$$\epsilon _{\cP (\cU)}:Ho(\cP(\cU))\to D(X).$$

\begin{lemma} $\epsilon _{\cP (\cU)}$ is an equivalence.
\end{lemma}

\begin{pf} Since $\cP(X)\hookrightarrow D(X)$ is an equivalence it follows that
$\epsilon _{\cP(\cU)} $
is essentially surjective. It remains to prove that $\epsilon _{\cP (\cU)}$ is full and faithful. Fix $P,Q
\in \cP(X)$. It suffices to prove that

1) $\Hom _{Ho(C(X))}(P,C_{\cU}(Q))=\Hom _{D(X)}(P,Q);$

2) $\Hom_{Ho(C(X))}(P,C_{\cU}(Q))=\Hom _{Ho(C(X))}(C_{\cU}(P),
C_{\cU}(Q)).$

 By devissage we may assume that complexes $P$ and $Q$ are vector bundles
places in degree 0. Let $i:U\hookrightarrow X$ be an embedding of an affine open subset.
Then for any $n$
$$\Hom _{D(X)}(P,i_*i^*Q[n])=\Hom _{D(U)}(i^*P,i^*Q[n])=0, \quad \text{if $n\neq 0$},$$
since $i^*P$ is a vector bundle on the affine variety $U$. This implies 1).

It remains to prove 2). Let $U,V\subset X$ be open subsets.
Note that $\Hom _{\cO_X}(P_U,Q_V)=0$ if $V\not\subset U$. Assume
that $V\subset U$. Then
$$\Hom _{\cO_X}(P_U,Q_V)=\Hom_{\cO_V}(P\vert _V,Q\vert _V)=
\Hom _{\cO_V}(\cO_X\vert _V,(P^*\otimes Q)\vert _V).$$
So we may assume that $P=\cO_X$. Let $Q_W$ be one of the summands  in the complex $C_{\cU}(Q)$.
(We assume that $Q_W$ is shifted to degree zero).
Then it follows that the complex $\Hom _{\cO_X}(C_{\cU}(\cO_X),Q_W)$ is acyclic except in degree zero
and
$$H^0(\Hom (C_{\cU}(\cO_X),Q_W))=\Gamma (W,Q)=\Hom(\cO_X,Q_W).$$
This proves the lemma.
\end{pf}

It follows from Lemma 5.4 that  $(\cP(\cU),\epsilon _{\cP (\cU)})$ is a standard enhancement
of $D(X)$.

Given smooth projective varieties $X$, $Y$ choose affine coverings $\cU$ and $\cV$ of $X$ and
$Y$ respectively. Then $\cU\times \cV$ is an affine covering of $X\times Y$. Given
$P\in \cP(X)$, $Q\in \cP(Y)$, $U\in \cU$, $V\in \cV$ we have $P_U\boxtimes Q_V=
(P\boxtimes Q)_{U\times V}$. This defines a DG functor
$$\boxtimes :\cP(\cU)\otimes \cP(\cV)\to \cP(\cU\times \cV).$$

\begin{lemma} The DG functor $\boxtimes $ is full and faithful.
\end{lemma}

\begin{pf} Let $P,P^\prime$ and $Q,Q^\prime$ be vector bundles on $X$ and $Y$ respectively
and fix
$U,U^\prime\in \cU$, $V,V^\prime\in \cV$. It suffices to prove that the natural map
$$\Hom (P_U,P^\prime_{U^\prime})\otimes \Hom (Q_V,Q^\prime _{V^\prime})\to
\Hom ((P\boxtimes Q)_{U\times V},(P^\prime\boxtimes Q^\prime)_{U^\prime\times V^\prime})$$
is an isomorphism, where all $\Hom $'s are taken in the usual categories of
quasi-coherent sheaves. Note that
both left and right hand sides are zero if $U^\prime \not\subset U$ or $V^\prime \not\subset
V$. So we may assume that $U^\prime \subset U$ and $V^\prime \subset
V$. Using the adjunction of direct and inverse image functors the question is reduced
to the following "affine" situation: Let $A$ and $B$ be noetherian $k$-algebras,
$M,M^\prime$ and
$N,N^\prime$ be
$A$- and $B$- modules respectively. Assume that the modules $M$ and $N$ are finitely
generated. Then the natural map
$$\Hom _A(M,M^\prime)\otimes _k\Hom _B(N,N^\prime)\to\Hom_{A\otimes B}(M\otimes N,M^\prime
\otimes N^\prime)$$
is an isomorphism. By taking resolutions
$$F_1\to F_0\to M\to 0,\quad \quad G_1\to G_0\to N\to 0,$$
where $F$'s and $G$'s are free modules of finite rank we may assume that $M=A$ and $N=B$,
in which case the assertion is clear.
\end{pf}

Since the $DG$ functor $\boxtimes $ is full and faithful the following diagram of DG functors
is commutative.
$$\begin{array}{ccc}
\cP(\cU)\otimes \cP(\cV) & \stackrel{h}{\longrightarrow} & (\cP(\cU)\otimes \cP(\cV))^0\text{-mod}\\
\boxtimes \downarrow & & \uparrow \Res _{\boxtimes}\\
\cP(\cU\times \cV) & \stackrel{h}{\longrightarrow} & \cP(\cU\times \cV)^0\text{-mod}
\end{array}
$$

The DG category $\cP(\cU\times \cV)$ is pretriangulated and perfect (since $D(X\times Y)$ is
Karoubian), so the DG functor $h:\cP(\cU\times \cV)  \to \cP(\cU\times \cV)^0\text{-mod}$
induces a
quasi-equivalence $\cP(\cU\times \cV)\to \text{Perf-}\cP(\cU\times \cV)$ (3.19).

Let $\cA \subset \cP(\cU\times \cV)$ be the smallest full pretriangulated DG subcategory, which
contains $\boxtimes (\cP(\cU)\otimes \cP(\cV))$. The commutativity of the above diagram
implies that the DG functor $\Res _{\boxtimes}\cdot h$ maps $\cA$ to
$\alpha((\cP(\cU)\otimes \cP(\cV))^{\text{pre-tr}})$ and the induces functor
$$Ho(\Res _{\boxtimes}\cdot h):Ho(\cA)\to Ho(\text{Perf-}(\cP(\cU)\otimes \cP(\cV)))$$
is
full and faithful.

By Theorems 1.9 and 1.10 the subcategory $\boxtimes (D(X)\times D(Y))$ classically (even strongly)
generates $D(X\times Y)$.
Therefore the triangulated subcategory $Ho(\cA)\subset Ho(\cP(\cU\times \cV))$
is
dense, and it follows that $\Res _{\boxtimes}\cdot h$ maps $\cP(\cU\times \cV)$ to
$\text{Perf-}(\cP(\cU)\otimes \cP(\cV))=\cP(\cU)\bullet \cP(\cV)$
and is a quasi-equivalence of these categories.
This proves the theorem.

The above proof of Theorem 5.5 gives us a more precise statement: standard enhancements are
compatible with products. Namely, we have the following proposition.

\begin{prop}  Let $X_1,...,X_n$ be smooth projective varieties. For each $i$ choose
a standard enhancement $(\cA _i,\epsilon _i)$ of $D(X_i)$.
Then there exists an equivalence of triangulated categories
$$\epsilon :Ho(\cA _1\bullet ...\bullet \cA _n)\to D(X_1\times ...\times X_n),$$
which makes
$(\cA _1\bullet ...\bullet\cA _n,\epsilon)$ a standard enhancement of
$D(X_1\times ...\times X_n)$ and
makes the following diagram commutative
$$\begin{array}{ccccc}
\times Ho(\cA _i) & \subset & Ho(\otimes \cA _i) & \stackrel{Ho(h)}{\longrightarrow} &
Ho(\cA _1\bullet ...\bullet \cA _n)\\
\downarrow \times \epsilon _i & & & & \downarrow \epsilon\\
\times D(X_i) &  & \stackrel{\boxtimes}{\longrightarrow} & &
D(X_1\times ...\times X_n).\\
\end{array}
$$
\end{prop}

 \begin{pf} Let $\cU _i$ be a finite affine covering of $X_i$. The last part of the proof
 of Theorem 5.5 implies the proposition for standard enhancements $\cA _i=\cP (\cU _i)$
 (the proof there is presented for $n=2$, but the general case is the same). The case of
 a general standard enhancement now follows from the functoriality of the product
 $\bullet $ in each variable (Remark 4.3).
 \end{pf}

\section{An application: representability of standard functors}

\begin{defn} Let $X$ and $Y$ be smooth projective varieties. A covariant (resp.
contravariant) functor  $F :D(X)\to D(Y)$
is {\it standard} if there exist standard enhancements $(\cA ,\epsilon _{\cA})$
and $(\cB ,\epsilon _{\cB})$ of $D(X)$ and $D(Y)$ respectively and a covariant (resp.
a contravariant) DG functor
$\tilde{F}:\cA \to \cB$ such that the functorial diagram
$$\begin{array}{ccc}
Ho(\cA) & \stackrel{Ho(\tilde{F})}{\longrightarrow} & Ho(\cB)\\
\epsilon _{\cA} \downarrow & & \downarrow \epsilon _{\cB}\\
D(X) & \stackrel{F}{\longrightarrow} & D(Y)
\end{array}
$$
is commutative up to an isomorphism.
We call $\tilde{F}$ a DG {\it lift} of $F$.
\end{defn}

Note that a standard functor is necessarily exact.

\begin{conj} All exact functors between $D(X)$ and $D(Y)$ are standard.
\end{conj}

\begin{example} The Serre functor $S=S_X:D(X)\to D(X)$ is standard. Indeed, the functor
$S_X$ is tensoring by a line bundle $\omega _X$ and then shifting by the dimension of $X$,
so it lifts for example to the standard enhancement $I(X)$.
\end{example}

Let $X$ be a smooth projective variety.
Recall the anti-involution ${\bf D}:D(X)\to D(X)$
$${\bf D}(\cdot )=\bbR \shHom _{\cO _X}(\cdot ,\cO _X).$$

\begin{lemma} The functor ${\bf D}$ is standard.
\end{lemma}

\begin{pf} Choose a finite affine covering $\cU$ of $X$ and consider the full DG subcategory
$\cP ^-(\cU)\subset C(X)$ defined similarly to $\cP(\cU)$ above except we use left
\v{C}ech resolutions instead of  right ones. Namely, given a perfect complex $T\in D(X)$
consider its resolution $C_{\cU}^-(T)$; it consists of direct sums of $\cO _X$-modules
of the form $i_!i^*T^j$, where $i:U\hookrightarrow X$ is the embedding of an open subset
$U$, which is the intersection of a few subsets from $\cU$. As in the case of $\cP (\cU)$
one proves that the tautological functor $\epsilon _{\cP^-(\cU)}:Ho(\cP^-(\cU))\to D(X)$ is
an equivalence, so that $(\cP^-(\cU),\epsilon _{\cP^-(\cU)})$ is a standard enhancement of
$D(X)$. Note that objects of $\cP^-(\cU)$ are acyclic for the functor $\shHom _{\cO _X}
(\cdot ,\cO_X)$.
Hence for a perfect complex $T$ we have
$${\bf D}(T)=\shHom(C_{\cU}^-(T),\cO_X).$$
But notice that the complex $\shHom(C_{\cU}^-(T),\cO_X)$ is equal to the complex
$C_{\cU}(\shHom(T,\cO _X))\in \cP(\cU)$. Thus we have lifted the duality ${\bf D}:D(X)\to
D(X)$ to a contravariant DG functor between the enhancements $\cP^-(\cU)$ and $\cP (\cU)$,
which shows that ${\bf D}$ is standard.
\end{pf}

\begin{lemma} The collection of standard covariant functors is closed under taking the (left and right)
adjoints, inverses (in case the functor is an equivalence) and taking composition of functors.
\end{lemma}

\begin{pf} Let $X$ and $Y$ be smooth projective varieties. Let $F:D(X)\to D(Y)$ be
an exact functor which is standard. More precisely, assume that there exist
standard enhancements $(\cA ,\epsilon _{\cA})$ and $(\cB ,\epsilon _{\cB})$ of
$D(X)$ and $D(Y)$ respectively and that there is  DG lift
$\tilde{F}:\cA \to \cB$ of the functor $F$.

Notice that the DG categories $\cA$ and $\cB$ are perfect, hence the Yoneda DG functors
$h:\cA \to \text{Perf-}\cA,$ $h:\cB \to \text{Perf-}\cB,$ are quasi-equivalences.
Thus $(\text{Perf-}\cA,\epsilon _{\cA}\cdot Ho(h)^{-1})$ and
$(\text{Perf-}\cB,\epsilon _{\cB}\cdot Ho(h)^{-1})$ are also standard enhancements of
$D(X)$ and $D(Y)$ respectively and
$$\Ind _{\tilde{F}}:\text{Perf-}\cA \to \text{Perf-}\cB$$
is another DG lift of $F$.
Its right adjoint DG functor $\Res _{\tilde{F}}$ is a DG lift of the right adjoint to
$F$. This proves that the right adjoint to $F$ and hence the inverse of $F$ (in case $F$
is an equivalence) are also standard.

Next we treat the composition of functors.
Namely, let $Z$ be another smooth projective variety  and
$G:D(Y)\to D(Z)$ be a standard functor with a DG lift $\tilde{G}: \cC \to \cD$ for
some standard enhancements $(\cC ,\epsilon _{\cC})$, $(\cD ,\epsilon _{\cD})$ of
$D(Y)$ and $D(Z)$ respectively. If there exists a quasi-equivalence of enhancements
$\mu :\cB \to \cC$, then the functor $G\cdot F$ has a DG lift $\tilde{G}\cdot \mu \cdot
\tilde{F}$. It there exists a quasi-equivalence of enhancements $\nu :\cC \to \cB$, then
the functor $G\cdot F$ has a DG lift
$$\Ind _{\tilde{G}}\cdot \Res _{\nu}\cdot \Ind _{\tilde{F}}:\text{Perf-}\cA\to
\text{Perf-}\cD.$$
Since in general there exists a chain of quasi-equivalences of standard enhancements connecting
$\cB$ and $\cC$ this shows that the functor $G\cdot F$ is standard.

It remains to show that a left adjoint to $F$ is standard. Notice that if $G$ is the
right adjoint to $F$ (which we proved is standard), then $S_X^{-1}\cdot G\cdot S_Y$ is
the left adjoint, which is therefore also standard.
\end{pf}

\begin{defn} Let $X$ and $Y$ be smooth projective varieties. An exact covariant functor
$F:D(X)\to D(Y)$ is {\it represented} by an object $P\in D(X\times Y)$ if there exists
an isomorphism of functors
$$F(\cdot )=\bbR q_*(p^*(\cdot)\stackrel{\bbL}{\otimes}P),$$
where $X\stackrel{p}{\leftarrow} X\times Y \stackrel{q}{\rightarrow} Y$ are the projections.
We say that $F$ is {\it representable} if it is represented by some object $P$.
\end{defn}

It is important to know that a given functor is representable, since, in particular,
it gives us an algebraic
cycle on $X\times Y$ (the Chern character applied to $P$), and construction of
algebraic cycles is always a difficult problem in algebraic geometry. Recall the
following important
theorem of Orlov [Or2].

\begin{thm} Let $X$ and $Y$ be smooth algebraic varieties and $F:D(X)\to D(Y)$ be
a covariant exact functor which is full and faithful. Then $F$ is representable.
\end{thm}

It is expected that the theorem holds without the assumption on $F$ being full and faithful.

\begin{thm} Let $X$ and $Y$ be smooth algebraic varieties and $F:D(X)\to D(Y)$ an
exact covariant functor which is standard. Then $F$ is representable.
\end{thm}

We need two lemmas.

\begin{lemma} In the notation of the theorem consider the (bi-contravariant)
functor $\theta :D(X)\times D(Y)\to Vect$ defined as follows:
$$\theta (A,B)=\Hom _{D(Y)}(B,F({\bf D}A)).$$
Assume that there exists a contravariant cohomological functor $\tau: D(X\times Y)\to Vect$
such that the functors $\theta $ and $\tau \cdot \boxtimes$ are isomorphic. Then $F$ is
representable.
\end{lemma}

\begin{pf} Since the category $D(X\times Y)$ is saturated there exists an object
$P\in D(X\times Y)$ and an isomorphism of functors $\tau (\cdot )=\Hom (\cdot ,P)$.
We have the following isomorphisms of functors ([Ha] II.5.15,5.14,5.11)
$$\begin{array}{rcl}
\Hom _{D(Y)}(B,F({\bf D}A)) & = & \theta (A,B)\\
& = & \Hom _{D(X\times Y)}(p^*A\stackrel{\bbL}{\otimes}q^*B ,P)\\
& = & \Hom _{D(X\times Y)}(q^*B, \bbR \shHom (p^*A ,P))\\
& = & \Hom _{D(X\times Y)}(q^*B, {\bf D}(p^*A)\stackrel{\bbL}{\otimes}P)\\
& = & \Hom _{D(Y)}(B, \bbR q_*({\bf D}(p^*A)\stackrel{\bbL}{\otimes}P))
\end{array}
$$
By [Ha] II.5.8 there is a functorial isomorphism ${\bf D}(p^*A)=
p^*({\bf D}(A))$. Summarizing we obtain the following functorial isomorphism
$$\Hom _{D(Y)}(B,F({\bf D}A))=
\Hom _{D(Y)}(B, \bbR q_*(p^*({\bf D}(A))\stackrel{\bbL}{\otimes}P)),$$
which implies an isomorphism of functors $F({\bf D}A)=
\bbR q_*(p^*({\bf D}(A))\stackrel{\bbL}{\otimes}P)$.
Since ${\bf D}^2=Id$ it follows that the functor $F$ is represented by the object
$P$.
\end{pf}

The next lemma explains the role of the assumption on the functor $F$ to be standard.

\begin{lemma} A standard functor $F$  in the above theorem satisfies the assumptions of
Lemma 6.9.
\end{lemma}

\begin{pf} Let $\tilde{F}:\cA \to \cB$ be a DG lift of $F$ where $(\cA,\epsilon _{\cA})$,
$(\cB ,\epsilon _{\cB})$ are standard enhancements of $D(X)$ and $D(Y)$ respectively.
We  know that the contravariant functor ${\bf D} :D(X)\to D(Y)$ has a DG lift
$\tilde{{\bf D}}:\cP^-(\cU)\to \cP(\cU)$ (Lemma 6.4). First let us show that
the composition $F\cdot {\bf D}$ is a standard functor. The argument is similar to the proof
of Lemma 6.5 above. Namely, if there exists a quasi-equivalence of enhancements
$\mu :\cP(\cU)\to \cA$, then
$$\tilde{F} \cdot \mu \cdot \tilde{{\bf D}}:\cP ^-(\cU)\to \cB$$
is a DG lift of $F\cdot {\bf D}$. If there exists a quasi-equivalence of enhancements
$\nu :\cA \to \cP(\cU)$, then
$$\Ind _{\tilde{F}}\cdot \Res _{\nu}\cdot h\cdot \tilde{{\bf D}}:\cP^-(\cU)\to \text{Perf-}\cB$$
is a DG lift of $F\cdot {\bf D}$. In general the standard enhancements $\cP(\cU)$ and $\cA$ of
$D(X)$ are connected by a chain of quasi-equivalences and we can use the above procedure
at each step to show that $G=F\cdot {\bf D}$ is standard. So we may assume that $G$
has a DG lift $\tilde{G}:\cC \to \cD$ for some standard enhancements
$(\cC ,\epsilon _{\cC})$, $(\cD ,\epsilon _{\cD})$ of $D(X)$ and $D(Y)$ respectively.

Consider the DG module $\tilde{\theta}\in (\cC \otimes \cD)^0\text{-mod}$ defined
by
$$\tilde{\theta}(A,B)=\Hom _{\cD}(B,\tilde{G}(A)).$$
Note that the functors $Ho(\tilde{\theta})$ and $\theta \cdot (\epsilon _{\cC}
\otimes \epsilon _{\cD})$ are isomorphic.

\medskip

\noindent{\it Claim.} There exists a DG module $\bar{\theta}\in (\cC \bullet \cD)^0\text{-mod}$
such that $\bar{\theta}\cdot h=\tilde{\theta}$, where $h:\cC \otimes \cD \to
\cC \bullet \cD$ is the
Yoneda embedding.

\medskip

Indeed, since $\cC \bullet \cD$ is a DG subcategory of $(\cC \otimes \cD)^0\text{-mod}$ we may
define
$$\bar{\theta}(Z)=\Hom_{(\cC \otimes \cD)^0\text{-mod}}(Z,\tilde{\theta}),\quad
\text{for $Z\in \cC\bullet \cD$}.$$
Notice that if $Z=h((A,B))$ for $(A,B)\in \cC \otimes \cD$, then
$$\Hom_{(\cC \otimes \cD)^0\text{-mod}}(Z,\tilde{\theta})=\tilde{\theta}((A,B)),$$
so we have $\bar{\theta}\cdot h=\tilde{\theta}$. This proves the claim.

By Proposition 5.9 there exists an equivalence $\epsilon:Ho(\cC\bullet \cD)\to D(X\times Y)$,
such that the diagram
$$\begin{array}{ccccc}
Ho(\cC)\times Ho(\cD) & \subset & Ho(\cC \otimes \cD) & \stackrel{Ho(h)}{\rightarrow} & Ho(\cC \bullet \cD)\\
\downarrow \epsilon _{\cC}\times \epsilon _{\cD} & & & & \downarrow \epsilon\\
D(X)\times D(Y) & & \stackrel{\boxtimes}{\rightarrow} & & D(X\times Y)
\end{array}
$$
is commutative. Therefore if we put $\tau =Ho(\bar{\theta})\cdot \epsilon ^{-1}$, then there
exists an isomorphism of functors $\tau \cdot \boxtimes =\theta$. This proves the lemma and
the theorem.
\end{pf}

\section{Another application: a motivic measure}

\subsection{The Grothendieck group $\Gamma$}

\begin{defn} Consider the abelian group $\Gamma =\Gamma _k$ generated by
quasi-equivalence classes
of pretriangulated categories $I(X)$ for smooth projective varieties $X$ over $k$ with
relations
coming from semiorthogonal decompositions. Namely, we impose the relation
$$[I(X)]=[I(Y_1)]+...+[I(Y_n)],$$
if there exist pretriangulated categories $\cA$, $\cB _1$,..., $\cB_n$,
which are quasi-equivalent to
$I(X)$, $I(Y_1)$,..., $I(Y_n)$ respectively and the following properties hold

1) $\cB_1$,..., $\cB _n$ are DG subcategories of $\cA$,

2) $Ho(\cB _1)$,..., $Ho(\cB _n )$ are admissible subcategories of $Ho(\cA)$,

3) $Ho(\cA)=(Ho(\cB _1),...,Ho(\cB _n))$ is a semiorthogonal decomposition.
\end{defn}

Notice that the group $\Gamma$ is defined similarly to the group $\cP\cT$, except we only
use quasi-equivalence classes of pretriangulated $I(X)$. In particular there is the
obvious group homomorphism $\beta :\Gamma \to \cP\cT$ (which may not be injective).
Just like in the case of $\cP\cT$ it follows (using Theorem 5.5) that $\Gamma $ is a
commutative associative ring. Note that $[I(pt)]=[DG(k)^f]$ is the identity in $\Gamma$.

\begin{remark} Semiorthogonal decompositions of derived categories $D(X)$ for smooth
projective $X$ tend to be rare. Examples are given in [BoOr] and in Propositions 7.3 and 7.5
below. Semiorthogonal summands of categories $D(X)$ are not always equivalent to
categories $D(Y)$ for a smooth projective $Y$. So it might
be more natural to include among generators of $\Gamma$ the  (quasi-equivalence classes of)
standard enhancements of semiorthogonal summands of categories $D(X)$.
\end{remark}

\begin{prop} Let $Z$ be a smooth projective variety and $E$ be a projectivization of a vector
bundle of dimension $d$ on $Z$. Then there is the following relation in $\Gamma$:
$$[I(E)]=d[I(Z)].$$
\end{prop}

\begin{pf} The proof of this and the next proposition is based on two theorems in
[Or].

Let $V\to Z$ be a vector bundle of dimension $d$ on $Z$ such that $p: E\to Z$ is its
projectivization. There exists a canonical relatively ample line bundle $\cO _V(1)$ on $E$.
Consider the functors $p_s^*:D(Z)\to D(E)$, where $p_s^*$ is the composition of the
inverse image functor $p^*=\bbL p^*:D(Z)\to D(E)$ with the tensoring by the line bundle
$\cO _V(s)$. The functor $p^*_0=p^*$ (and hence also the functors $p^*_s$
for all $s$) is full and faithful [Or]. Denote by $D(Z)_s\subset D(E)$ the strictly full
subcategory of $D(E)$ which is the essential image of the functor $p^*_s$. It is proved in [Or]
that $(D(Z)_{-d+1},...,D(Z)_0)$ is a semiorthogonal decomposition of $D(E)$.

The functor $p^*$ induces a fully faithful DG functor $p^*=p^*_0 :I(Z)\to C(E)$.
 Consider the
commutative diagram of functors
$$\begin{array}{ccc}
Ho(I(Z)) & \stackrel{Ho(p^*)}{\longrightarrow} & Ho(p^* (I(Z)))\\
\downarrow & & \downarrow \\
D(Z) & \stackrel{p^*}{\longrightarrow} & D(E),
\end{array}
$$
where the vertical arrows are the obvious functors. This shows that the right vertical
arrow is full and faithful, because the other three are such. Hence $p^* (I(Z))$
(which is DG equivalent to $I(Z)$)
is a standard enhancement of the category $D(Z)_0\subset D(E)$.

For each $s$ we may consider the DG functor $p^* _s :I(Z)\to C(E)$ which is the
composition of $p^* $ with tensoring by the line bundle $\cO _V(s)$. Clearly, the
image DG category $p^* _s(I(Z))$ (which is DG equivalent to $I(Z)$)
is a standard enhancement of the subcategory
$D(Z)_s\subset D(E)$. The proposition follows.
\end{pf}

\begin{example} Taking $Z=pt$ we deduce that $[I(\bbP ^{d-1})]=d$ in $\Gamma$. This, of course,
 also follows from the theorem of Beilinson on the resolution of the diagonal in
$\bbP ^n\times \bbP ^n$.
\end{example}

\begin{prop} Let $Y$ be a smooth projective variety and $Z\subset  Y$ a smooth
closed subvariety of codimension $d$. Denote by $X$ the blowup of $Y$ along $Z$ with
the exceptional divisor $E\subset X$. Then there is the following relation in the group
$\Gamma$:
$$[I(X)]+[I(Z)]=[I(Y)]+[I(E)].$$
\end{prop}

\begin{pf} We have the following obvious commutative diagram
$$\begin{array}{rcl}
E & \stackrel{j}{\hookrightarrow} & X\\
p\downarrow & &\downarrow \pi \\
Z & \hookrightarrow & Y,
\end{array}
$$
where $E$ is the projectivization of a vector bundle of dimension d over $Z$. As in the proof
of the last lemma consider the full subcategories $D(Z)_{-d+1},...,D(Z)_0\subset D(E)$.
Denote by $j_{*s}$ the restriction of the direct image functor $j_*:D(E)\to D(X)$ to
the subcategory $D(Z)_s$. The following statements are proved in [Or]:

1) the functors $j_{*s}$ are full are faithful for all $s$;

2) the inverse image functor $\pi ^*:D(Y)\to D(X)$ is full and faithful;

3) $(j_{*(-d+1)}D(Z)_{-d+1},..., j_{*(-1)}D(Z)_{-1},\pi^* D(Y))$
is a semiorthogonal decomposition
of $D(X)$.

Hence (as in the proof of the last proposition) it follows that the  DG subcategory
$j_{*s}p^*_s(I(Z))\subset C(X)$ (which is DG equivalent to $I(Z)$)
is a standard enhancement of the subcategory
$j_{*s}D(Z)_s\subset D(X)$. Also $\pi ^*(I(Y))$ (which is DG equivalent to $I(Y)$)
is a standard enhancement of $\pi ^*(D(Y))$.
We obtain the following relation in $\Gamma$
$$[I(X)]=[I(Y)]+(d-1)[I(Z)].$$
It remains to apply the last proposition.
\end{pf}

\subsection{A motivic measure}

Let $\cV_k$ denote the collection of $k$-varieties. The Grothendieck group $K_0[\cV _k]$ is
generated by isomorphism classes of $k$-varieties with relations
$$[X]=[Y]+[X-Y],$$
if $Y\subset X$ is a closed subvariety.
The product of varieties over $k$ makes $K_0[\cV _k]$ a commutative associative ring.
Needless to say that it is interesting and important to understand the structure of this
ring. However, very little is known. For example, only recently it was proved [Po] that
 $K_0[\cV _k]$ has zero divisors for any field $k$.
The following result gives a description of the quotient ring $K_0[\cV _{\bbC}]/
\langle[\bbL]\rangle$, where $\bbL$ is the class of the affine line $\bbA ^1$.
Consider the multiplicative monoid $SB$ of stable birational equivalence classes varieties
over $\bbC$. Let $\bbZ[SB]$ denote the corresponding monoid ring.

\begin{thm} [LaLu] There exists a natural isomorphism of rings
$$ K_0[\cV _{\bbC}]/
\langle[\bbL]\rangle\simeq \bbZ[SB].$$
\end{thm}

A ring homomorphism $K_0[\cV _k]\to A$ is called a {\it motivic measure}. We claim
that there exists a natural (surjective) motivic measure
$$\mu:K_0[\cV _{\bbC}]\to \Gamma _{\bbC},$$
whose kernel contains the principal ideal generated by $\bbL-1$ (1=[pt]).

Indeed, recall the Loojenga-Bittner presentation [Lo],[Bit] of the ring $K_0[\cV _{\bbC}]$:
it is generated by isomorphism classes of smooth projective varieties with the defining set of
relations
$$[X]+[Z]=[Y]+[E],$$
if $X$ is the blowup of $Y$ along $Z$ with the exceptional divisor $E$.
Then by Proposition 7.5 and Theorem 5.5 the correspondence $[X]\mapsto [I(X)]$
(for smooth peojective $X$) extends to a well defined ring homomorphism
$\mu:K_0[\cV _{\bbC}]\to \Gamma _{\bbC}.$ By Example 7.4 $\mu ([\bbP ^1])=2$, hence
$\mu (\bbL)=1$.

The induced ring homomorphism
$$\bar{\mu}:K_0[\cV _{\bbC}]/\langle\bbL -1\rangle \to \Gamma _{\bbC}$$
is probably not injective (for example, by Mukai's theorem [Mu] the categories
$D(A)$ and $D(\hat{A})$ are equivalent if $A$ is an abelian variety and $\hat{A}$ is its dual),
but seems to be close to injective.

\end{document}